\documentclass[11pt]{amsart}
\usepackage{amssymb,mathrsfs}
\usepackage[author-year]{amsrefs}

\title{A Coupling, and the Darling--Erd\H{o}s Conjectures}
   \thanks{The research of D. Kh.\ was
   supported by a grant from the National Science Foundation}
   \author[D. Khoshnevisan]{Davar Khoshnevisan}
   \address{Department of Mathematics, The University of Utah,
      155 S. 1400 E. Salt Lake City, UT 84112--0090}
   \email{davar@math.utah.edu}
   \urladdr{http://www.math.utah.edu/\~{}davar}
   \author[D. A. Levin]{David A. Levin}
   \address{Department\@ of Mathematics, The University of Oregon,
      Eugene, OR 97403--1221}
   \email{dlevin@uoregon.edu}
   \urladdr{http://www.uoregon.edu/\~{}dlevin}
\keywords{The Darling--Erd\H{o}s theorem, coupling, strong approximations}
\subjclass[2000]{60G70, 60F15}
\date{August 16, 2005}
\date{September 2, 2005}
\date{October 6, 2005}
\date{November 14, 2005}
\date{November 16, 2005}

\theoremstyle{plain}{
\newtheorem{theorem}{Theorem}[section]}
\theoremstyle{plain}{
\newtheorem{proposition}[theorem]{Proposition}}
\theoremstyle{plain}{
   \newtheorem{lemma}[theorem]{Lemma}}
\theoremstyle{plain}{
   }
\theoremstyle{definition}{
   }
\theoremstyle{definition}{
   }
\theoremstyle{remark}{
   \newtheorem{remark}[theorem]{Remark}}
\theoremstyle{definition}{
   }
\theoremstyle{plain}{
\newtheorem*{DE}{The Darling--Erd\H{os} Conjecture}}
\usepackage{amssymb,euscript}
\usepackage[author-year]{amsrefs}
\numberwithin{equation}{section}

\newcommand{\e}{\epsilon}
\renewcommand{\P}{\mathrm{P}}
\newcommand{\E}{\mathrm{E}}
\newcommand{\T}{U}
\newcommand{\1}{{\bf 1}}
\begin{document}
\begin{abstract}
	We derive a new coupling of the running maximum of
	an Ornstein--Uhlenbeck process and the
	running maximum of an explicit i.i.d.\ sequence.
	We use this coupling to verify a conjecture of Darling and Erd\H{o}s
	\ycite{DE}.
\end{abstract}
\maketitle

\section{Introduction and Main Results}

Let $\{\xi_i\}_{i=1}^\infty$ be a sequence of independent,
identically distributed random variables 
with $\E[\xi_1]=0$ and $\E[\xi_1^2]=1$, and define
\begin{equation}
	\T_n := \max_{1\le k\le n} 
	\left(\frac{\xi_1 + \cdots + \xi_k}{k^{1/2}}\right)
	\qquad(n\ge 1).
\end{equation}
According to the law of the iterated logarithm (LIL)
of \ocite{hartmanwintner},
\begin{equation}
	\limsup_{n\to\infty} \frac{\T_n}{(2\log \log  n)^{1/2}}=1
	\quad\text{almost surely.}
\end{equation}
Here and throughout, ``$\log x$'' and/or
``$\log(x)$'' act as short-hand for ``$\ln(x\vee e)$.''
In a remarkable paper \ycite{DE}, Darling and Erd\H{os}
establish the following variation of the LIL.
 \begin{theorem}\label{thm:DE}
	 Assume that $\E[\xi_1]=0$, $\E[\xi_1^2]=1$,
	 and $\E\{|\xi_1|^3\}<\infty$. Then for all real
	 numbers $x$,
	 \begin{equation}\begin{split}
		&\lim_{n\to\infty} \P\left\{
			a(n) \T_n - b(n) \le x \right\} =\exp
			\left(- \frac{e^{-x}}{(4\pi)^{1/2}}
			\right),\quad\text{where}\\
		&a(x) := (2\log\log x)^{1/2},\text{ and }
			b(x) := 2\log\log x +\frac12\log\log\log x.
	 \end{split}\end{equation}
 \end{theorem}
Subsequently, \ocite{oodaira} and \ocite{shorack}
improved the integrability condition to
``$\E\{|\xi_1|^{2+\e}\}<\infty$ for some $\e>0$.''
The definitive result, along these lines, is due
to \ocite{einmahl89} who improved the integrability
condition further to ``$\E\{\xi_1^2\log \log |\xi_1|\}<\infty$,'' and proved that
the said condition is optimal. Related works can be found
in \ocite{bertoin} and \ocite{einmahlmason}.

At the very end of their paper, Darling and Erd\H{o}s pose
two conjectures about strong-limit analogues of
Theorem \ref{thm:DE}. Define
\begin{equation}\begin{split}
	c_1 & := \limsup_{n\to\infty}\frac{a(n)
		\T_n - b(n)}{\log\log\log\log n},\\
	c_2 & := -\liminf_{n\to\infty}\frac{a(n)
		\T_n - b(n)}{\log\log\log\log n}.
\end{split}\end{equation}
By the Kolmogorov zero-one law, $c_1$ and $c_2$
are constants almost surely.

\begin{DE}\label{DE1}
	With probability one, $0<c_1,c_2<\infty$.
\end{DE}

The main goal of this paper is to prove that
$c_1=\infty$ and $c_2=1$, under a mild
moment condition on the distribution of
$\xi_1$.
Thus, half of the conjecture is true while the other half
is false.

Our proof involves first deriving
a new and novel coupling (Theorem \ref{coupling})
of the running maximum of an Ornstein--Uhlenbeck
process with the running maximum of a certain i.i.d.\ process.
Theorem \ref{thm:DE} follows readily from this coupling.
Our solution to the Darling--Erd\H{o}s conjecture also follows from it,
but requires a little more work.
En route, we present also an integral test  (Theorem \ref{thm:test})
for the lower envelope of
the Ornstein--Uhlenbeck process; see also our integral
test for pure-cosine lacunary series (Theorem
\ref{thm:lac}). We adapt Breiman's
terminology \ycite{breiman}, and assert that
our integral tests are ``very delicate.''
See Remark \ref{delicate} below for an explanation.

Next we  present the precise form 
of our solution to the Darling--Erd\H{o}s conjecture.

\begin{theorem}\label{rw1}
	With probability one,
	$c_1=\infty$ and $c_2=1$, as long as
	\begin{equation}\label{eq:integrability}
		\E[\xi_1]=0,\
		\E[\xi_1^2]=1,\ \text{and}\
		\E\left[ \xi_1^2\log\log|\xi_1| \right]<\infty.
	\end{equation}
\end{theorem}

It turns out that the seemingly more interesting contradictory
half (i.e., $c_1=\infty$) is, in fact, not very deep. It relies
only on results that were known prior to the work of
Darling and Erd\H{o}s \citelist{\cite{erdos}\cite{feller}},
and does not require any further computations.
On the other hand, developing the formula ``$c_2=1$'' seems to 
require some new ideas.

\ocite{shorack} has observed that the strong invariance
principles of \ocite{PhilippStout} yield
Darling--Erd\H{o}s theorems for many processes with dependent
increments as well. Next, we explore this observation further in the special
case of lacunary pure-cosine series only. It is both possible 
and tempting to use
the constructions of \ocite{PhilippStout}*{Sections 3--12} and
\ocite{berkes}, and find further
embellishments (lacunary series with weights),
and other applications (functions of strongly mixing random variables,
partial sums of stationary Gaussian processes with
long-range dependence, Markov sequences, etc.).
However, we will not do that because no 
further ideas are needed
to carry out such a program. Thus, we complete the
Introduction by stating some implications of
the present work in the context of pure-cosine
series of the lacunary type.

Let $\{n_k\}_{k=1}^\infty$ be a sequence of real numbers that
satisfies the \emph{lacunarity condition},
\begin{equation}\label{eq:lac}
	\liminf_{k\to\infty} \frac{n_{k+1}}{n_k} >1.
\end{equation}
Consider the pure-cosine series,
\begin{equation}
	f_n(\omega)= 
	2^{1/2}\sum_{1\le k\le n} \cos(2\pi n_k\omega),
\end{equation}
where $n\ge 1$ and $0\le \omega\le 1$.
\ocite{shorack} has proved that the
Darling--Erd\H{o}s Theorem \ref{thm:DE} continues to hold
if we replace $\T_n$ by 
$\max_{1\le k\le n}(f_k/k^{1/2})$, and 
the measure $\P$ by
the Lebesgue measure on $[0\,,1]$. More precisely,
Shorack's theorem asserts that for all real numbers $x$, 
the following holds: 
\begin{equation}\label{eq:shorack}\begin{split}
	&\lim_{n\to\infty} \left| \left\{ 0\le \omega\le 1:\ 
		a(n)\max_{1\le k\le n}\left( \frac{f_k(\omega)}{k^{1/2}} 
		\right) -b(n)\le x \right\}\right|\\
	&= \exp\left( -\frac{e^{-x}}{(4\pi)^{1/2}} \right),
\end{split}\end{equation}
where $|\cdots|$ denotes the Lebesgue measure.
Shorack's proof rests on a strong approximation theorem of
\ocite{PhilippStout}*{Theorem 3.1, p.\ 12} and a verification of
two asymptotic negligibility conditions
\cite{shorack}*{eq.'s (2.1) and (2.3)}. 
It is possible, and not too hard, to replace the said asymptotic
negligibility conditions with an appeal to the Erd\H{o}s--Feller
integral test, and use the strong approximation theorem
of \ocite{berkes} to prove the following.

\begin{proposition} 
	Equation \eqref{eq:shorack} continues to hold if the lacunarity condition
	\eqref{eq:lac} is
	replaced by the following weaker hypothesis:
	\begin{equation}\label{eq:berkes}
		\text{There exists $0\le \alpha<\tfrac12$ such that \ 
		$\lim_{k\to\infty}
		k^\alpha \left( \frac{n_{k+1}}{n_k} - 1
		\right) =\infty$.}
	\end{equation}
\end{proposition}
[Mimic the argument the begins with \eqref{r0} below.]

We obtain a corresponding integral test as
an immediate consequence of the analysis
of the present paper.
 Before we describe the said integral test
 let us define
\begin{equation}
	\mathscr{F}_n(\omega) := a(n) \max_{1\le k\le n}
	\left( \frac{f_k(\omega)}{k^{1/2}} \right) 
	- b(n) + \frac{\log(4\pi)}{2},
\end{equation}
where $n= 1,2,\ldots$ and $0\le\omega\le 1$.
The following is an immediate consequence
of the strong approximation theorem of \ocite{berkes}
and our Theorem \ref{thm:test} below.

\begin{theorem}\label{thm:lac}
	Suppose $\{n_k\}_{k=1}^\infty$ satisfies \eqref{eq:berkes}
	and $h:[0\,,\infty)\to[0\,,\infty)$ is non-decreasing
	and satisfies 
	$h(n+1)-h(n)=o(1/n)$ as $n\to\infty$.
	Then,
	\begin{equation}
		\left|  \left\{ 0\le \omega\le 1:\
		\mathscr{F}_n(\omega) \le - h(n) \text{ i.o.} \right\}\right| 
		= 
		\begin{cases}
			0&\text{ if $\mathbf{J}(h)=\infty$,}\\
			1&\text{ if $\mathbf{J}(h)<\infty$.}
		\end{cases}
	\end{equation}
	Here, ``i.o.'' means ``infinitely often,'' and
	\begin{equation}\label{eq:h}
		\mathbf{J}(h) := \int_1^\infty \exp\left(
		h(t) - e^{h(t)} \right) \frac{dt}{t}.
	\end{equation}
\end{theorem}

\noindent\textbf{Acknowledgement.} 
This paper owes much of its existence to Zhan Shi, and
Rodrigo Ba\~{n}uelos generously
pointed out to us the connections between our coupling
theorem and 
processes with dependent increments such as lacunary series.
We heartily thank them both.

\section{A Coupling for OU Processes}

Throughout, $X:=\{X_t\}_{t\ge 0}$ denotes the 
Ornstein--Uhlenbeck (OU) process. We recall that $X$ is a continuous, centered,
Gaussian process with $\textrm{Cov}( X_s\,,X_t ) = \exp(-|t-s|/2)$
for $s,t\ge 0$.
We recall also that $X$ is a stationary and ergodic diffusion, and
$X_0$ is standard normal.

Let $\ell:=\{\ell_t\}_{t\ge 0}$ denote the local times of $X$ at zero.
It is well known that $\ell$ is continuous (a.s.), and 
\begin{equation}\label{eq:LT}
	\ell_t := \lim_{\e\to 0} \frac1\e \int_0^t \1_{
	\{ 0< X_s< \e\}}\, ds
	\qquad(t\ge 0).
\end{equation}
The convergence holds almost surely and in 
$L^p(\P)$ for all $p\in[1\,,\infty)$.
Define
\begin{equation}
	\tau(t) := \inf\{ s>0:\ \ell_s>t\}
	\qquad(t\ge 0).
\end{equation}
Also, we introduce the process
$M:=\{M_n\}_{n=1}^\infty$ as follows:
\begin{equation}
	M_n := \sup_{\tau(n-1)\le s\le \tau(n)} X_s
	\qquad\text{for all }n\ge1.
\end{equation}
By the strong Markov property, $M$ is an i.i.d.\
sequence. Also,
the distribution of $M$ has been computed explicitly
in Proposition 2.2 of Khoshnevisan, Levin, and Shi
\ycite{KLS}. It reads as follows:
For all real numbers $x$,
\begin{equation}\label{eq:cdf}
	\P\{ M_1 \le x\} := F(x) := \exp\left( -\frac{1}{
	2\int_0^{\max(x,0)} \exp(y^2/2)\, dy}\right),
\end{equation}
where $\exp(-1/0):=0$.
Now we present and prove the following coupling.

\begin{theorem}\label{coupling} 
	As $t\to\infty$,
	\begin{equation}
		\P\left\{ \sup_{0\le s\le t} X_s \neq
		\max_{1\le j\le (2\pi)^{-1/2}t} M_j
		\right\} =O\left( \left( 
		\log  (t) /t\right)^{1/2}\right).
	\end{equation}
\end{theorem}

Our proof of Theorem \ref{coupling} requires a technical
lemma. But first, let us observe from \eqref{eq:LT} that
$\E[\ell(t)] = t(2\pi)^{-1/2}$.
Thus, by the ergodic theorem,
$\lim_{t\to\infty} \ell(t)/t = (2\pi)^{-1/2}$ a.s.
A time substitution then yields the following:
\begin{equation}\label{tau/t}
	\lim_{t\to\infty} \frac{\tau(t)}{t} = (2\pi)^{1/2}
	\quad\text{a.s.}
\end{equation}
Convergence holds also in $L^p(\P)$ for all $p\in[1\,,\infty)$;
confer with \eqref{tau:MGF} below.
The aforementioned technical lemma is the following
quantitative refinement of the ergodic theorem \eqref{tau/t}.

\begin{proposition}\label{pr:MD}
	For all $\alpha>0$ there exists $\beta>1$
	such that for all $t\ge 1$,
	\begin{equation}
		\P\left\{ \sup_{0\le s\le t}
		\left|\tau(s) - s(2\pi)^{1/2}\right|
		\ge\beta (t\log  t)^{1/2} \right\} \le
		\beta t^{-\alpha}.
	\end{equation}
\end{proposition}

\begin{proof}
	The functions $\tau$, $f(x)=x$,  and
	$g(x)=(x \log  x)^{1/2}$ are non-decreasing.
	Therefore, it suffices to prove that for all
	$\alpha>0$ there exists $\beta>1$ such that
	for all integers $n\ge 1$,
	\begin{equation}\label{goal:MD}
		\P\left\{ \max_{1\le j\le n}
		\left|\tau(i) - i(2\pi)^{1/2}\right| 
		\ge \beta (n\log   n)^{1/2}\right\} \le
		\beta n^{-\alpha}.
	\end{equation}
	By the strong Markov property, $\tau$ is
	an ordinary random walk. This, \eqref{tau/t},
	and the Kolmogorov strong law of large numbers
	together imply that $\E[\tau(1)]=(2\pi)^{1/2}$.
	We propose to verify that $\tau(1)$ has a finite
	moment generating function. 
	Then, \eqref{goal:MD} follows at once
	from the classical moderate deviations estimates
	of Cram\'er \cite{il}*{Theorem 6.1.1, p.\ 156}.
	
	Define $T_1:=\inf\{ s>0:\ X_s=1\}$, $T_1':=
	\inf\{s>T_1:\ X_s=0\}$, $T_2:=\inf\{
	s>T_1':\ X_s=1\}$, etc. These are the respective
	crossing-times of one and zero.
	
	Because the speed measure
	of $X$ decays faster than exponentially, $T_1$
	has exponential moments of all order; i.e.,
	$\E_x[\exp(aT_1)]<\infty$ for all $a>0$ and 
	$x\in(-\infty\,,\infty)$
	\cite{Mandl}*{Lemma 2, p.\ 112}. Furthermore,
	if $T_0:=0$, then the strong Markov property of $X$ guarantees
	that:
	\begin{enumerate}
	\item
		$\{(T_n-T_{n-1}\,,\ell(T_n)-\ell(T_{n-1}))\}_{n=2}^\infty$
		is an independent sequence;
	\item
		The $\P$-distribution of $T_n-T_{n-1}$ is the same as
		the $\P_1$-distribution of $T_1$;
	\item
		The $\P$-distribution of $\ell(T_n)-\ell(T_{n-1})$,
		and the $\P_1$-distribution of $\ell(T_1)$, are 
		the same, in fact exponential.
	\end{enumerate}
	Define $\mu:=\E_1[\ell(T_1)]$, and note the bounds,
	\begin{equation}\begin{split}
		\P\{ \tau(1) > T_n \} & =
			\P\left\{\ell(T_n) < 1 \right\}\\
		&\le \P\left\{ \sum_{2\le i\le n}
			\left( \ell(T_i) -\ell(T_{i-1})\right)
			<1 \right\}\\
		&\le \P\left\{ \sum_{2\le i\le n}
			\left( \ell(T_i) -\ell(T_{i-1}) \right)\le
			\frac{n-1}{2} \mu \right\},
	\end{split}\end{equation}
	for all $n$ large enough. By large deviations, this implies
	the existence of a constant $c$ such that for all $n$ large,
	$\P\{ \tau(1) > T_n \} \le \exp(-cn)$; see, for example
	\ocite{il}*{eq.\ 13.2.1, p.\ 245}.
	Thus, for all $n$ large, letting $\nu := \E_1[T_1]$,
	\begin{equation}\begin{split}
		\P\left\{ \tau(1) > 2\nu n\right\} &
			\le \exp(-cn) + \P\left\{ T_n\ge 2\nu n\right\}.
	\end{split}\end{equation}
	We can write $T_n=(T_n-T_{n-1})+\cdots +(T_2-T_1)+T_1$.
	We recall that $T_1$ and $T_i-T_{i-1}$ are independent,
	and the $(T_i-T_{i-1})$'s are i.i.d.\ Moreover, all have finite 
	exponential moments of all orders. Therefore,
	another appeal to large deviations proves that
	$\P\{T_n>2\nu n\}\le \exp(-c'n)$ for some constant
	$c'>0$ that does not depend on $n$. From this we
	can conclude that 
	\begin{equation}\label{tau:MGF}
		\E\left[\exp(a\tau(1))\right]<\infty
		\quad\text{for some }a>0.
	\end{equation}
	As was mentioned earlier, \eqref{goal:MD} follows
	from this at once. This proves the proposition.
\end{proof}

Now we can prove Theorem \ref{coupling}

\begin{proof} [Proof of Theorem \ref{coupling}]
	Define $E_\beta(t)$ to be complement of the event,
	\begin{equation}
		\left\{  \max_{1\le j\le t -\beta(t\log  t)^{1/2}} M_j  \le
		\sup_{0\le s\le t(2\pi)^{1/2}} X_s
		\le \max_{1\le j\le t +\beta(t\log  t)^{1/2}} M_j  \right\}.
	\end{equation}
	For all integers $n,m\ge 1$,
	\begin{equation}
		\P\left\{ \max_{1\le j\le n}M_j
		\neq \max_{1\le i\le n+m}M_i \right\} = \frac{m}{n+m}.
	\end{equation}
	Therefore, as $t\to\infty$,
	\begin{equation}\begin{split}
		&\P\left\{ \sup_{0\le s\le t(2\pi)^{1/2}}
			X_s \neq \max_{1\le j\le t}M_j\right\}\\
		&\quad \le \P\left( E_\beta(t) \right)
			+ \P\left\{ \max_{1\le j\le t -\beta(t\log  t)^{1/2}} M_j
			\neq \max_{1\le j\le t +\beta(t\log  t)^{1/2}} M_j  \right\}\\
		&\quad = \P\left( E_\beta(t) \right)
			+ \left( 2 \beta+o(1) \right)
			\left( \frac{\log  t}{t}\right)^{1/2}.
	\end{split}\end{equation}
	Thanks to Proposition \ref{pr:MD}, we can choose
	and fix $\beta$, once and for all, so large that the probability
	of $E_\beta(t)$ is
	$O((\log (t)/t)^{1/2})$. This completes our proof. 
\end{proof}

\section{An Integral Test for OU Processes}

By \eqref{eq:cdf} and a direct computation,
\begin{equation}
	\lim_{n\to\infty} \P\left\{ \max_{1\le j\le n}
	M_j  \le \left( 2\log  n+\log \log  n +x \right)^{1/2}
	\right\} = \exp\left( - \frac{e^{-{x/2}}}{2^{1/2}}
	\right),
\end{equation}
for every real number $x$.
Our coupling (Theorem \ref{coupling}) then yields
the following without further effort: For all real numbers $x$,
 \begin{equation}\label{eq:DE}
	 \lim_{t\to\infty} \P\left\{
	 a\left(e^t\right) \sup_{0\le s\le t}X_s - b\left(e^t\right)
	  \le x \right\}
	 =\exp\left( - \frac{e^{-x}}{(4\pi)^{1/2}}\right).
\end{equation}
This is the analogue of the Darling--Erd\H{o}s
theorem for OU processes, and is implicitly the
first part of the original proof of Theorem \ref{thm:DE}.
Explicitly, it appears as a special case of a result of
Pickands \ycite{pickands}*{Theorem 4.4} for
Gaussian processes.  Earlier, 
Newell \ycite{newell}*{pp. 491\ndash 492}, studying diffusions,
obtained a version of \eqref{eq:DE} asymptotic in $x$.
It also can be found in \ocite{shorack}, and as a consequence
of a much more general theorem of \ocite{bertoin}*{Theorem 3}.
There is extensive literature on the maximum of stationary Gaussian
processes; for a sampling see \ocite{vr}, \ocite{cramer},
\ocite{qw}, and \ocite{berman}.


The main purpose of this section is to derive an integral
test that corresponds to the $\liminf$ behavior of 
$a\left(e^t\right)\sup_{0\le s\le t}X_s-b\left(e^t\right)$.
  We find it more convenient
to work with the following variant:
\begin{equation}
	\mathscr{X}(t):=a(e^t) \sup_{0\le s\le t}X_s - b(e^t)
	 +\frac{\log (4\pi)}{2}\qquad(t>0).
\end{equation}

Suppose $g$ is a Borel-measurable function that
is non-decreasing ultimately, and has
the following additional properties:
\begin{equation}\label{eq:phi}
	\frac{\log (g(n+1))}{\log (g(n))} = 1+ o(1/n)
	\qquad(n\to\infty).
\end{equation}
Also define for all measurable functions $g:[1\,,\infty)\to(0\,,\infty)$,
\begin{equation}
	\mathbf{I}(g) := \int_1^\infty \frac{\log (g(t))}{tg(t)}\, dt.
\end{equation}
To compare with \eqref{eq:h} we note merely that
$\mathbf{I}(g)=\mathbf{J}(\log\log g)$.
Then we have the following:

\begin{theorem}\label{thm:test}
	Assume $g$ satisfies \eqref{eq:phi}, and define
	$F$ to be the event that the random set
	$\{t>0:\ \mathscr{X}(t)\le - \log \log  (g(t)) \}$ is unbounded.
	Then $\P(F)=0$ or $1$ according as $\mathbf{I}(g)=\infty$
	or $\mathbf{I}(g)<\infty$.
\end{theorem}

\begin{remark}\label{delicate}
	This is a ``very delicate'' LIL \cite{breiman} in the following sense:
	If we perturb the gauge function $-\log  \log (g(t))$ even a little
	and replace it by $-\log ( c+\log (g(t)))$, then the end-result could
	be vastly different. For example, it follows from Theorem
	\ref{thm:test} that for all $\e>0$, a.s.:
	\begin{equation}\label{LLL}\begin{split}
		\mathscr{X}(t)& \le -\log \left( \log \log  t + 2
			\log \log \log  t \right) \text{ unboundedly, whereas}\\
		\mathscr{X}(t)& >-\log \left( \log \log  t + (2+\e)\log \log \log  t
			\right) \text{ eventually}.
	\end{split}\end{equation}
	The delicateness of the integral test is now seen, for the
	difference between the right-most terms in \eqref{LLL}
	is $(\e+o(1)) (\log \log \log  t)(\log \log  t)^{-1}=o(1)$ as $t\to\infty$,
	whereas $\lim_{t\to\infty}\mathscr{X}(t)=-\infty$ a.s.
\end{remark}

We will prove that Theorem \ref{thm:test} is a consequence
of two technical lemmas. Those are developed first.
Throughout, $g$ is a Borel-measurable function that
is non-decreasing ultimately.

\begin{lemma}\label{step1}
	If $g$ satisfies \eqref{eq:phi}, then
	$\mathbf{I}(g) <\infty$ if and only if
	almost surely for all but a finite number of
	$n$'s,
	\begin{equation}
		\max_{1\le j\le n} M_j
		> \left( 2\log  n + \log \log  n -\log  (2) - 2\log \log (g(n))
		\right)^{1/2}.
	\end{equation}
\end{lemma}

\begin{proof}
	We can assume, without loss of generality, that
	\begin{equation}\label{eq:gg}
		\log  n \le g(n) \le \log  n \cdot (\log \log  n)^3\qquad
		(n\ge 1).
	\end{equation}
	For otherwise we could replace $g$ everywhere by $g_1$,
	where
	\begin{equation}
		g_1(x):=\min\left\{ 
		\max( g(x),\log  x) \,,\,\log  x\cdot(\log \log  x)^3 \right\}.
	\end{equation}
	
	Recall $F$ from \eqref{eq:cdf}, and define
	$\bar{F}:=1-F$.
	Let $\{u_n\}_{n=1}^\infty$ be a non-decreasing sequence
	of positive real numbers such that $n\bar{F}(u_n)$ is
	also non-decreasing. Then,
	according to Theorem 4.3.1 of \ocite{galambos}*{p.\ 214},
	\begin{equation}
		\P\left\{ \max_{1\le j\le n}M_j < u_n \text{ i.o.} \right\}=
		\begin{cases}
		0&\text{if }\sum_{n= 1}^\infty \bar{F}(u_n)
			e^{-n\bar{F}(u_n)} <\infty,\\
		1&\text{if }\sum_{n= 1}^\infty \bar{F}(u_n) 
			e^{-n\bar{F}(u_n)}  =\infty.
		\end{cases}
	\end{equation}
	In fact, the monotonicity of $u_n$ and
	$n\bar{F}(u_n)$ can be replaced
	by the following condition, as can be seen by inspecting the
	proofs in \ocite{galambos}*{pp.\ 214--222}:
	\begin{equation}\label{reduce}
		u_k\text{ is ultimately increasing, and }
		\max_{1\le k\le n} k \bar{F}(u_k)
		= n\bar{F}(u_n) + O(1),
	\end{equation}
	as $n\to\infty$.
	A little bit of calculus shows that
	\begin{equation}
		\left| \int_0^x \exp(y^2/2)\, dy -
		\frac 1x e^{x^2/2} \right| = O\left(
		x^{-3} e^{-x^2/2} \right)\quad(x\to\infty).
	\end{equation}
	Therefore, by Taylor's expansion,
	\begin{equation}\label{barF}
		\left| \bar{F}(x) - \frac12 x e^{-x^2/2}\right| = 
		O\left( x^{-3} e^{-x^2/2} \right)
		\quad(x\to\infty).
	\end{equation}
	We apply the preceding with
	\begin{equation}
		u_n := \left( 2\log  n + \log \log  n -\log  (2) - 2\log \log  (g(n))
		\right)^{1/2} \quad(n\ge 1).
	\end{equation}
	Then, \eqref{eq:phi} and \eqref{eq:gg} together imply that
	\begin{equation}
		\sum_{n\ge 1} \bar{F}(u_n) 
		e^{-n\bar{F}(u_n)} <\infty\ \Longleftrightarrow\
		\mathbf{I}(g)<\infty.
	\end{equation}
	We omit the details as they involve routine computations.
	Similar work shows that we have \eqref{reduce}. 
	Whence follows the lemma.
\end{proof}

We apply Theorem \ref{coupling} and Lemma
\ref{step1} in conjunction to obtain the following.

\begin{lemma}\label{step2}
	Under the preceding conditions, 
	$\mathbf{I}(g)<\infty$ if and only if
	almost surely for all but a bounded set of $t$'s,
	\begin{equation}
		\sup_{0\le s\le t}X_s
		> \left( 2\log  t + \log \log  t -\log (4\pi) - 2\log \log (g(t))
		\right)^{1/2}.
	\end{equation}
\end{lemma}

\begin{proof}
	Let $m_n:=[\exp(\tau n/\log  n)]$ where $\tau\in(0\,,1)$
	is small enough but otherwise fixed; see the top of page 222 of
	\ocite{galambos} for details.
	
	Thanks to \ocite{galambos}*{pp.\ 218--222},
	our proof of Lemma \ref{step1} implies, in fact,
	that  $\mathbf{I}(g)<\infty$ if and only if
	almost surely for all but a finite number of $n$'s,
	\begin{equation}\begin{split}
		&\max_{1\le j\le m_n}M_j\\
		&> \left( 2\log  m_n + \log \log  m_n -\log  (2) - 
			2\log \log (g(m_n)) \right)^{1/2},
	\end{split}\end{equation}
	and this is, in turn, equivalent to the validity of the following
	for all but a finite number of $n$'s:
	\begin{equation}\begin{split}
		&\max_{1\le j\le m_{n-1}}M_j\\
		&> \left( 2\log  m_n + \log \log  m_n -\log  (2) - 
			2\log \log (g(m_n)) \right)^{1/2},
	\end{split}\end{equation}
	But according to Theorem \ref{coupling} and the Borel--Cantelli
	lemma, with probability one,
	\begin{equation}
		\max_{1\le j\le m_n}M_j  \ =
		\sup_{0 \le s\le m_n(2\pi)^{1/2}} X_s
		\quad\text{eventually as $n\to\infty$}. 
	\end{equation}
	This, and monotonicity, together prove the lemma.
\end{proof}

\begin{proof}[Proof of Theorem \ref{thm:test}]
	As in the proof of Lemma \ref{step1}, we assume
	without loss of generality that \eqref{eq:gg} holds.
	Next we make some real-variable computations.
	
	For all $\e>0$ small enough,
	\begin{equation}\label{ep}
		(1+\e)^{1/2} \le 1+\frac \e2 \le 
		(1+\e+\e^2)^{1/2}.
	\end{equation}
	Choose and fix a real number $p$ and define
	\begin{equation}
		\e(t) := \frac{\log \log  t}{2\log  t} +
		\frac{p}{2\log  t} - \frac{\log \log (g(t))}{\log t}.
	\end{equation}
	Plug in $\e:=\e(t)$ in the first bound in \eqref{ep} 
	to deduce that
	\begin{equation}\label{bd1}\begin{split}
		&\left( 1 + \frac{\log \log  t}{2\log  t} +
			\frac{p}{2\log  t} - \frac{\log \log (g(t))}{\log  t} \right)^{1/2}\\
		&\le  1 + \frac{\log \log  t}{4 \log  t}
			+ \frac{p}{4\log  t} - \frac{\log \log (g(t))}{2\log  t}.
	\end{split}\end{equation}
	Condition \eqref{eq:gg} implies that
	\begin{equation}
		\e^2 (t) = O\left( \left( \frac{\log \log  t}{\log t} \right)^2 \right)
		=o\left( \frac{1}{\log  t 
		\cdot \log (g(t))}\right) \quad(t\to\infty).
	\end{equation}
	Now choose and fix a number $c\in(0\,,1)$, and apply
	\eqref{eq:phi} once again to find that
	\begin{equation}
		0\le \frac{\log \log (g(t)) - \log \log (cg(t))}{\log  t} 
		=  \frac{\log (1/c)+o(1) }{\log  t\cdot \log (g(t))}
		\quad(t\to\infty).
	\end{equation}
	The last two displays together with \eqref{ep} imply that 
	almost surely as $t$ grows to infinity,
	\begin{equation}\label{bd2}
		\left( \frac{\log \log  t}{2\log  t} +
		\frac{p}{2\log  t} - \frac{\log \log (cg(t))}{\log  t} \right)^{1/2}
		\ge 1 + \frac{\e(t)}{2}
		\quad\text{eventually.}
	\end{equation}
	The theorem follows readily from
	\eqref{bd1}, \eqref{bd2}, and Lemma \ref{step2},
	because $\mathbf{I}(g)<\infty$ if and only if
	$\mathbf{I}(cg)<\infty$.
\end{proof}

\section{Proof of Theorem \ref{rw1}}

We can recast Feller's \ycite{feller} improvement on the integral test 
of \ocite{erdos} as follows:
For any non-decreasing sequence $\{r_n\}_{n=1}^\infty$
of positive real numbers,
\begin{equation}\label{eq:feller}
	\P\left\{ \T_n \le r_n 
	\ \text{eventually}\right\} =1 \quad
	\text{if and only if}\quad
	\sum_{n\ge 1} \frac{ r_n}{n}e^{- r_n^2/2}<\infty.
\end{equation}
Feller's test is valid solely
under the condition \eqref{eq:integrability};
see \ocite{einmahl89} where a fatal gap in
Feller's proof was bridged. Of course, if the probability in \eqref{eq:feller}
is strictly less than one, then it is zero. 
This follows from the Kolmogorov zero-one law.

We can apply Feller's test to the sequence
\begin{equation}
	 r_n := \left( 2\log\log n +
	3\log\log\log n + \theta \log\log\log\log n\right)^{1/2}.
\end{equation}
For this particular choice of $ r_n$'s, 
the summability condition of \eqref{eq:feller} holds if and only if
$\theta>2$. A little algebra yields \eqref{c_1:corrected}, 
whence follows that $c_1=\infty$.
Now we turn to proving that $c_2=1$, all the time assuming that
\eqref{eq:integrability} holds.

First of all, we note that 
\begin{equation}\label{r0}
	\liminf_{t\to\infty} \frac{a(t)\sup_{0\le s\le \log t}X_s
	- b(t)}{\log\log\log\log t} =-1\quad\text{almost surely.}
\end{equation}
This is a consequence of Theorem \ref{thm:test},
and was mentioned earlier in a slightly different
form (Remark \ref{delicate}). From here on, we use strong approximations.
The forthcoming argument is inspired by those of
\ocite{oodaira} and \ocite{shorack}.

Define $S_0:=\xi_0:=0$
and $S_t:=\sum_{0\le j\le t}
\xi_j$ as usual ($t> 0$).
We can note that $H(t):=t^2\log\log  t$
satisfies the conditions
of Theorem 2 of \ocite{einmahl}. Thus,  
Einmahl's theorem implies that we can construct
$\{S_t\}_{t\ge 0}$ together with a (standard)
Brownian motion $\{B_t\}_{t\ge 0}$---on
a suitably-chosen probability space---such that
almost surely,
\begin{equation}\label{SA}
	|S_s- B_s | =o\left( \left(\frac{s}{
	\log\log s} \right)^{1/2}\right)
	\qquad(s\to\infty).
\end{equation}
Define
\begin{equation}
	c(n) := \exp\left( \frac{\log n}{(\log\log n)^3} \right)
	\qquad(n\ge 1).
\end{equation}

Choose and fix $\e>0$.
According to \eqref{eq:feller}, the following holds
almost surely: For all $n$ sufficiently large,
\begin{equation}\begin{split}
	\T_{c(n)} &\le \left( 2\log\log c(n) + 3\log\log\log n
		+(2+\e)\log\log\log\log n
		\right)^{1/2}\\
	&\le a(n) \left[ 1 - \frac{3\log\log\log n}{2\log\log n}
		+ \frac{(1+\e)\log\log\log\log n}{\log\log n}\right].
\end{split}\end{equation}
Consequently,
\begin{equation}\label{r1}
	\liminf_{n\to\infty} \frac{a(n)\T_{c(n)}-b(n)}{
	\log\log\log\log n} =-\infty\qquad\text{almost surely.}
\end{equation}
Similarly, we have
\begin{equation}\label{r2}
	\liminf_{n\to\infty} \frac{a(n)\sup_{1\le s\le c(n)}(B_s/
	s^{1/2})-b(n)}{
	\log\log\log\log n} =-\infty\qquad\text{almost surely.}
\end{equation}
On the other hand, thanks to \eqref{SA}, almost surely,
\begin{equation}
	\sup_{c(n)\le s\le n} \frac{|B_s-S_s|}{s^{1/2}}
	= o\left( \left( \log\log c(n) \right)^{-1/2}\right).
\end{equation}
Because $\log\log c(n)=(1+o(1))\log\log n$, 
it follows that almost surely,
\begin{equation}\label{r3}
	\limsup_{t\to\infty} \frac{a(n)}{\log\log\log\log n}
	\sup_{c(n)\le s\le n} \frac{|S_s-B_s|}{s^{1/2}}=0.
\end{equation}
But the process $\{e^{-s/2}B_{\exp(s)}\}_{s\ge 0}$ has
the same finite-dimensional distributions as $X$.
Combine this observation with \eqref{r0}, \eqref{r1}, \eqref{r2},
and \eqref{r3} to complete the proof.

\section{Epilogue}
\subsection{Further Refinements}
During the course of our proof of Theorem \ref{rw1}, we 
have proved the following facts:
\begin{enumerate}
	\item
		If $\E\{ \xi_1^2 \log \log |\xi_1| \}<\infty$, then
		with probability one,
		\begin{equation} \label{c_1:corrected}
			\limsup_{n \rightarrow \infty} \frac{(2\log\log
			n)^{1/2} \T_n -
			2\log\log n - \frac{3}{2} \log\log\log n}{ 
			\log\log\log\log n} = 1. 
		\end{equation}
	\item
		The moment condition \eqref{eq:integrability}
		implies also that with probability one,
		\begin{equation}\label{Eq:liminf}
			\liminf_{n \rightarrow \infty}
			\frac{(2\log\log n)^{1/2}\T_n
			- 2\log\log n - \tfrac12\log\log\log n}{
			\log\log\log\log n} = -1 . 
		\end{equation}
\end{enumerate}
It is not hard to see from our arguments
that more stringent moment conditions
yield more detailed results. We leave the details
to the interested reader.

\subsection{Toward a Conjecture of Pickands}
There are abundant techniques already available
for handling the limit superior behavior of
the running maxima of Gaussian and diffusion processes
\cites{berman,pickands,pickands:lil,pickands:apmsgp,qw:abgp,qw}.
In comparison, the
literature on limit inferior behavior is scant.
The following noteworthy result along these lines is due to
\ocite{pickands:apmsgp}*{Theorem 3.2}:
\begin{equation} \label{Eq:PickandsLiminf}
	\liminf_{t \rightarrow \infty}
	\frac{ (2\log t)^{1/2}\sup_{0\le s\le t}X_s - 
	2\log t}{\log\log t}
	\geq \frac{1}{2} \quad\text{almost surely}.
\end{equation}
Our integral test (Theorem \ref{thm:test})
improves this in a definitive manner. In particular, 
it follows that with probability one,
\begin{equation} 
	\liminf_{t \rightarrow \infty}
	\frac{ (2\log t)^{1/2}\sup_{0\le s\le t}X_s - 
	2\log t}{\log\log t}
	= \frac{1}{2}.
\end{equation}
This verifies a conjecture of \ocite{pickands:apmsgp}*{p.\ 86}
in the special case of OU processes.
The  methods of \ocite{pickands:apmsgp} 
have been applied to study similar
problems in other settings
\cite{pickands:lil}. Our techniques are quite different,
however, and do not seem to work when the process
in question is not Markovian.


\begin{bibdiv}
\begin{biblist}

\bib{berkes}{article}{
    author={Berkes, I.},
     title={An almost sure invariance principle for lacunary trigonometric
            series},
   journal={Acta Math. Acad. Sci. Hungar.},
    volume={26},
      date={1975},
     pages={209\ndash 220},
}

\bib{berman}{book}{
    author={Berman, Simeon M.},
     title={Sojourns and Extremes of Stochastic Processes},
 publisher={Wadsworth \& Brooks/Cole},
     place={Pacific Grove, CA},
      date={1992},
}

\bib{bertoin}{article}{
    author={Bertoin, Jean},
     title={Darling-Erd\H os theorems for normalized sums of i.i.d.
            variables close to a stable law},
   journal={Ann. Probab.},
    volume={26{\it (2)}},
      date={1998},
     pages={832\ndash 852},
}

\bib{breiman}{article}{
    author={Breiman, Leo},
     title={A delicate law of the iterated logarithm for non-decreasing
            stable processes},
   journal={Ann. Math. Statist.},
    volume={39},
      date={1968},
     pages={1818\ndash 1824},
}
\bib{cramer}{article}{
  author = {Cram{\'{e}}r, H},
  title = {A limit theorem for the maximum value of certain stochastic
  processes},
  journal = {Th. Probab. Appl.},
  volume = {10},
  pages = {126\ndash 128},
  year = {1965},
}

\bib{DE}{article}{
    author={Darling, D. A.},
    author={Erd\H{o}s, P.},
     title={A limit theorem for the maximum of normalized sums of
            independent random variables},
   journal={Duke Math. J.},
    volume={23},
      date={1956},
     pages={143\ndash 155},
}

\bib{einmahl89}{article}{
    author={Einmahl, Uwe},
     title={The Darling-Erd\H os theorem for sums of i.i.d.\ random
            variables},
   journal={Probab.\ Theory Related Fields},
    volume={82{\it (2)}},
      date={1989},
     pages={241\ndash 257},
}


\bib{einmahl}{article}{
    author={Einmahl, Uwe},
     title={Strong invariance principles for partial sums of independent
            random vectors},
   journal={Ann. Probab.},
    volume={15{\it (4)}},
      date={1987},
    pages={1419\ndash 1440},
}

\bib{einmahlmason}{article}{
    author={Einmahl, Uwe},
    author={Mason, David M.},
     title={Darling-Erd\H os theorems for martingales},
   journal={J. Theoret. Probab.},
    volume={2{\it (4)}},
      date={1989},
     pages={437\ndash 460},
}

\bib{erdos}{article}{
  author =   {Erd\H{o}s, Paul},
  title =    {On the law of the iterated logarithm},
  journal =  {Ann.\ Math.},
  year =     {1942},
  volume =   {43{\it (2)}},
  pages =    {419--436},
}

\bib{feller}{article}{
    author={Feller, W.},
     title={The law of the iterated logarithm for identically distributed
            random variables},
   journal={Ann. of Math. (2)},
    volume={47},
      date={1946},
     pages={631\ndash 638},
}

\bib{galambos}{book}{
    author={Galambos, Janos},
     title={The Asymptotic Theory of Extreme Order Statistics}, 
   publisher={John Wiley \& Sons, New York-Chichester-Brisbane},
      date={1978},
}

\bib{hartmanwintner}{article}{
    author={Hartman, Philip},
    author={Wintner, Aurel},
     title={On the law of the iterated logarithm},
   journal={Amer.\ J. Math.},
    volume={63},
      date={1941},
     pages={169\ndash 176},
}


\bib{il}{book}{
    author={Ibragimov, I. A.},
    author={Linnik, Yu. V.},
     title={Independent and Stationary Sequences of Random Variables},
      note={Translation from the Russian edited by J. F. C. Kingman},
 publisher={Wolters-Noordhoff Publishing, Groningen},
      date={1971},
}

\bib{KLS}{article}{
     author =   {Khoshnevisan, Davar},
     author =   {Levin, David A.},
     author =   {Shi, Zhan},
      title =   {An extreme-value analysis of the LIL
		      for Brownian motion},
       year =   {2005},
    journal = {Electr.\ Comm.\ in Probab.},
     volume = {10, {\rm paper 20}},
       year = {2005},
      pages = {196--206},
}
%
\bib{Mandl}{book}{
    author={Mandl, Petr},
     title={Analytical Treatment of One-Dimensional Markov Processes},
 publisher={Academia Publishing House of the Czechoslovak Academy of
            Sciences, Prague},
      date={1968},
}

\bib{newell}{article}{
  author = {Newell, G.F.},
  title = {Asymptotic extreme value distributions for one dimensional
  diffusions}, 
  journal = {J. Math. Mech.},
  volume = {11},
  pages = {481\ndash 496},
  year = {1962},
}

\bib{oodaira}{incollection}{
    author={Oodaira, Hiroshi},
     title={Some limit theorems for the maximum of normalized sums of weakly
            dependent random variables},
 booktitle={Proceedings of the Third Japan-USSR Symposium on Probability
            Theory (Tashkent, 1975)},
     pages={1\ndash 13. Lecture Notes in Math., Vol. 550},
 publisher={Springer},
     place={Berlin},
      date={1976},
}

\bib{PhilippStout}{article}{
    author={Philipp, Walter},
    author={Stout, William},
     title={Almost {S}ure {I}nvariance {P}rinciples for {P}artial {S}ums of {W}eakly
            {D}ependent {R}andom {V}ariables},
   journal={Mem. Amer. Math. Soc. 2},
      date={1975},
    volume={161},
}

\bib{pickands}{article}{
  author = {Pickands, James, III},
  year = {1967},
  title = {Maxima of stationary Gaussian processes},
  journal = {Z. Wahrscheinlichkeitstheorie verw. Geb.},
  volume = {7},
  pages = {190\ndash 233},
}

\bib{pickands:lil}{article}{
  author = {Pickands, James, III},
  year = {1969a},
  title = {An Iterated Logarithm Law for the Maximum in a Stationary
  Gaussian Sequence},
  journal = {Z. Wahrscheinlichkeitstheorie verw. Geb.},
  volume = {12},
  pages = {344\ndash 353},
}
\bib{pickands:apmsgp}{article}{
  author = {Pickands, James, III},
  year = {1969b},
  title = {Asymptotic properties of the maximum in a stationary
  Gaussian process},
  journal = {Trans. Amer. Math. Soc.},
  pages = {75\ndash 86},
  volume = {145},
}

\bib{qw:abgp}{article}{
  author = {Qualls, Clifford},
  author = {Watanabe, Hisao},
  title = {An asymptotic 0-1 behavior of Gaussian processes},
  journal = {Ann. Math. Stat.},
  volume = {42},
  number = {6},
  year = {1971},
  pages = {2029\ndash 2035},
}

\bib{qw}{article}{
  author = {Qualls, Clifford},
  author = {Watanabe, Hisao},
  title = {Asymptotic Properties of Gaussian Processes},
  journal = {Ann. Math. Stat.},
  volume = {43},
  number = {2},
  year = {1972},
  pages = {580\ndash 596},
}

\bib{vr}{article}{
  author = {Volkonskii, V.A.},
  author = {Rozanov, Y.A},
  title = {Some limit theorems for random functions II},
  journal = {Th. Probab. Appl.},
  volume = {6},
  pages = {186\ndash 198},
  year = {1961},
}

\bib{shorack}{article}{
    author={Shorack, Galen R.},
     title={Extension of the Darling and Erd\H os theorem on the maximum of
            normalized sums},
   journal={Ann. Probab.},
    volume={7{\it (6)}},
      date={1979},
     pages={1092\ndash 1096},
}

%
\end{biblist}
\end{bibdiv}
\end{document}